\documentclass[12pt]{article}
\usepackage{amssymb,amsmath}
\def\hang{\hangindent\parindent}
\def\tex#1{\indent\llap{[#1]\enspace}\ignorespaces}
\def\re{\par\hang\tex}
\def\a{\alpha}
\def\b{\beta}
\def\d{\delta}
\def\D{\Delta}
\def\UU{{\mathcal U}}
\def\WW{{\mathcal W}}
\def\LL{{\mathcal L}}
\def\SS{{\mathcal{L}}(\G)}
\def\g{\gamma}
\def\G{\Gamma}
\def\Der{{\rm Der}}
\def\Inn{{\rm Inn}}

\def\Ker{{\rm Ker}}
\def\es{\varepsilon}
\def\Im{{\rm Im}}

\def\v{\varphi}

\def\ssc{\scriptscriptstyle}

\def\cl{\centerline}

\def\rar{\rightarrow}

\def\bs{\backslash}

\def\vs{\vspace*}

\def\ra{\rangle}

\def\la{\langle}
\def\ni{\noindent}
\def\ptl{\partial}

\def\Z{\mathbb{Z}{\ssc\,}}

\def\C{\mathbb{C}{\ssc\,}}
\def\F{\C}
\def\QED{\hfill$\Box$}
\numberwithin{equation}{section}
\newtheorem{theo}{Theorem}[section]
\newtheorem{defi}[theo]{Definition}

\newtheorem{lemm}[theo]{Lemma}
\newtheorem{prop}[theo]{Proposition}
\newtheorem{clai}{Claim}
\newtheorem{subc}{Subclaim}
\def\adddot{$\!\!\!${\bf.}\ \ }
\def\addbra{$\!\!\!${\bf)}\ \ }

\begin{document}
\cl{{\Large \bf
 Lie bialgebras of generalized Virasoro-like type}\footnote
{Supported by a NSF grant 10471096 of China, ``One Hundred Talents
Program'' from University of Science and Technology of China and
``Trans-Century Training Programme Foundation for the Talents'' from
National Education Ministry of China.}} \vs{6pt}

\cl{Yuezhu Wu$^{\,1,2)}$,  Guang'ai Song$^{\,3)}$, Yucai
Su$^{\,4)}$}

 \cl{\small
$^{1)}$Department of Mathematics, Shanghai Jiao Tong University,
 Shanghai 200240, China}

\cl{\small $^{2)}$Department of Mathematics, Qufu Normal University,
 Qufu 273165, China}

\cl{\small $^{3)}$College of Mathematics and Information
\vs{-3pt}Science, Shandong Institute of Business and} \cl{\small
Technology, Yantai, Shandong 264005, China}

 \cl{\small \small $^{4)}$Department of
Mathematics, University of Science and Technology of \vs{-3pt}China}
\cl{\small Hefei 230026, China}

 \cl{\small E-mail:  ycsu@ustc.edu.cn} \vs{6pt}

{\small
\parskip .005 truein
\baselineskip 3pt \lineskip 3pt

\noindent{{\bf Abstract.} In two recent papers by %Song and Su, and
the authors, all Lie bialgebra structures on %graded and nongraded
Lie algebras of generalized Witt type are classified. In this paper
all Lie bialgebra structures on generalized Virasoro-like algebras
are determined. It is proved that all such Lie bialgebras are
triangular coboundary. \vs{5pt}

}

}
\parskip .001 truein\baselineskip 8pt \lineskip 8pt

\vs{6pt}

\cl{\bf\S1. \
Introduction}\setcounter{section}{1}\setcounter{equation}{0} The
notion of Lie bialgebras was first introduced by Drinfeld in 1983
[D1] (cf.~[D2]) in a connection with quantum groups. Since then
there appeared a number of papers on Lie bialgebras (e.g., [M1, M2,
NT, N,
%S,
SS, WS, T]). %, which arise naturally in the study of Hamiltonian
%mechanics and Poisson Lie group.
Michaelis [M1] presented a class of %infinite-dimensional
Lie
bialgebras containing the Virasoro algebra (this type of Lie
bialgebras was classified by Ng and Taft [NT], cf.~[N, T]) and gave
a method on how to obtain the structure of a triangular coboundary
Lie bialgebra on a Lie algebra containing two elements $a,b$
satisfying $[a,b]=b$.

In two recent papers [SS, WS], %(see also [S])
all Lie bialgebra structures on Lie algebras of generalized Witt
type are classified. In this paper we shall determine all Lie
bialgebra structures on a class of Lie algebras
(cf.~(\ref{gvir-equa-gvir-like})), referred to as the {\it
generalized Virasoro-like algebras} (the structure and
representation theories of the Virasoro-like algebra have attracted
some authors' attentions because of its close relation with the
Virasoro algebra, e.g., [LT, MJ, X2, X3, ZM, ZZ]).

Let us recall the definition of Lie bialgebras. For a vector space
$\LL $ over the complex field $\F$,
%a field $\F$ of characteristic zero,
we define the {\it twist map} $\tau$ of $\LL  \otimes \LL $ and the
{\it cyclic map} $\xi$ of $\LL  \otimes \LL  \otimes \LL $
\vs{-6pt}by
\begin{equation}
\tau:x\otimes y\mapsto y \otimes x, \ \ \ \ \ \ \xi: x \otimes
y\otimes z\mapsto y \otimes z \otimes x \mbox{ \ \ \ for \ \
}x,y,z\in \LL \vs{-6pt}.
\end{equation}
Then a {\it Lie algebra} can be defined as a pair $(\LL ,\v)$
consisting of a vector space $\LL $ and a bilinear map $\v :\LL
\otimes \LL  \rar \LL $ (the {\it bracket} of $\LL $) satisfying the
following conditions\vs{-6pt},
\begin{eqnarray}
\label{gvir-Lie-s-s} \!\!\!\!\!\!\!\!\!\!\!\!&&
\Ker(1-\tau) \subset \Ker\,\v \mbox{ \ (skew-symmetry),}\\
\label{gvir-Lie-j-i} \!\!\!\!\!\!\!\!\!\!\!\!&& \v \cdot (1 \otimes
\v ) \cdot (1 + \xi +\xi^{2}) =0 : \ \LL  \otimes \LL  \otimes \LL
\rar \LL \mbox{ \ (Jacobi identity),}
\end{eqnarray}
where $1$ is the identity map of $\LL  \otimes \LL $. A {\it Lie
coalgebra} is a pair $(\LL , \D)$ consisting of a vector space $\LL
$ and a linear map $\D: \LL  \to \LL  \otimes \LL $ ({\it cobracket}
of $\LL $) satisfying the following conditions:
\begin{eqnarray}
\label{gvir-cLie-s-s} \!\!\!\!\!\!\!\!\!\!\!\!&&
\Im\,\D \subset \Im(1- \tau) \mbox{ \ (anti-commutativity),}\\
\label{gvir-cLie-j-i} \!\!\!\!\!\!\!\!\!\!\!\!&& (1 + \xi +\xi^{2})
\cdot (1 \otimes \D) \cdot \D =0:\ \LL  \to \LL  \otimes \LL
\otimes \LL \mbox{ \ (Jacobi identity).}
\end{eqnarray}
\begin{defi}\adddot
\rm A {\it Lie bialgebra} is a triple $(\LL ,\v, \D )$ such that
$(\LL ,\v)$ is a Lie algebra and $(\LL ,\D)$ is a Lie coalgebra and
the following {\it compatibility condition} holds:
\begin{eqnarray}
\label{gvir-bLie-d} \!\!\!\!\!\!\!\!\!\!\!\!&& \mbox{$\D  \v (x, y)
= x \cdot \D y - y \cdot \D x$ \ \ \ for \ \ $x, y \in \LL $, }
\end{eqnarray}
where the symbol ``$\cdot$'' means the action
\begin{equation}
\label{gvir-diag} x\cdot (\mbox{$\sum\limits_{i}$} {a_{i} \otimes
b_{i}}) = \mbox{$\sum\limits_{i}$} ( {[x, a_{i}] \otimes b_{i} +
a_{i} \otimes [x, b_{i}]})
\end{equation}
for $x, a_{i}, b_{i} \in \LL $, and in general $[x,y]=\v(x,y)$ for
$x,y \in \LL $.
\end{defi}

One shall notice that the significant difference between Lie
bialgebras and (associative) bialgebras lies in the compatibility
condition (\ref{gvir-bLie-d}): A bialgebra requires that
 $\D$ is an algebra
morphism: $ \D \cdot \v = (\v \otimes \v) \cdot (1 \otimes \tau
\otimes 1) \cdot \D \otimes \D, $ while a Lie bialgebra requires
that $\D$ is a derivation (cf.~(\ref{gvir-deriv})) of $\LL  \rar \LL
\otimes \LL .$ Thus the properties of Lie bialgebras are not similar
to those of bialgebras.
\begin{defi}\adddot
\label{gvir-def2} \rm (1) A {\it coboundary Lie bialgebra} is a
$(\LL , \v, \D,r),$ where $(\LL , \v, \D)$ is a Lie bialgebra and $r
\in \Im(1 - \tau) \subset \LL  \otimes \LL $ such that $\D$ is a
{\it coboundary of $r$}, i.e. $\D=\D_r$, where in general $\D_r$
(which is an inner derivation, cf.~(\ref{gvir-inn})) is defined by,
\begin{equation}
\label{gvir-D-r} \D_r (x) = x \cdot r \mbox{ \ \ for \ \ }x \in \LL
.
\end{equation}

(2) A coboundary Lie bialgebra $(\LL , \v,\D, r)$ is {\it
triangular} if it satisfies the following {\it classical Yang-Baxter
Equation} (CYBE):
\begin{equation}
\label{gvir-CYBE} c(r)=0,
\end{equation}
 where $c(r)$ is defined by
\begin{equation}
\label{gvir-add1-} \mbox{$c(r) = [r^{12} , r^{13}] +[r^{12} ,
r^{23}] +[r^{13} , r^{23}],$}
\end{equation}
and $r^{ij}$ are defined as follows: Denote $\UU(\LL )$ the
universal enveloping algebra of $\LL $ and $1$ the identity element
of $\UU (\LL )$. If $r =\sum_{i} {a_{i} \otimes b_{i}} \in \LL
\otimes \LL $, then
$$\begin{array}{l}
 r^{12} = r\otimes 1=\sum \limits_{i}{a_{i} \otimes b_{i}
\otimes 1,} \\[12pt]
r^{13}= (1\otimes \tau)(\tau\otimes 1)=\sum \limits_{i} {a_{i}
\otimes 1 \otimes b_{i}}, \\[12pt]
r^{23} =1\otimes r= \sum \limits_{i}{1 \otimes a_{i} \otimes
b_{i}},
\end{array}$$
are all elements in $\UU (\LL ) \otimes \UU (\LL ) \otimes \UU(\LL
).$

\end{defi}

Let us state our main results below. For any {\it nondegenerate}
additive subgroup $\G$ of $\F^2$ (namely, $\G$ contains a $\F$-basis
of $\F^2$), the {\it generalized Virasoro-like algebra} $\SS$ is a
Lie algebra with basis $\{L_\a,\ptl_1,\ptl_2\,|\,\a\in\G\bs\{0\}\}$
and bracket
\begin{equation}\label{gvir-equa-gvir-like}
[L_\a,L_\b]=(\a_1\b_2-\b_1\a_2)L_{\a+\b},\ \ \
[\ptl_i,L_\a]=\a_iL_\a \mbox{ \ for \ }\a,\b\in\G,\ i=1,2,
\end{equation}
where we use the convention that if an undefined notation appears in
an expression, we always treat it as zero; for instance, $L_\a=0$ if
$\a=0$. In particular, when $\G=\Z^2$, the derived subalgebra
$[\LL(\Z^2),\LL(\Z^2)]={\rm span}\{L_\a\,|\,\a\in\Z^2\bs\{0\}\}$ is
the (centerless) Virasoro-like algebra (e.g., [LT, MJ, ZZ]). The Lie
algebra $\SS$ is closely related to the Lie algebras of Block type
(cf.~[DZ, X3, Z1]) and the Lie algebras of Cartan type $S$ (cf.~[SX,
X1, Z2]).

For a Lie algebra $\LL$ and an $\LL$-module $V$, denote by
$H^1(\LL,V)$ the {\it first cohomology group} of $\LL$ with
coefficients in $V$. It is well-known that
\begin{equation}
H^1(\LL,V)\cong\Der(\LL,V)/\Inn(\LL,V),
\end{equation}
where $\Der(\LL,V)$ is the set of \textit{derivations} $d:\LL\to V$
which are linear maps satisfying
\begin{equation}
\label{gvir-deriv} d([x,y])=x\cdot d(y)-y\cdot d(x)\mbox{ \ for \
}x,y\in \LL,
\end{equation}
and $\Inn(\LL,V)$ is the set of {\it inner derivations} $a_{\rm
inn},\,a\in V$, defined by
\begin{equation}
\label{gvir-inn} a_{\rm inn}:x\mapsto x\cdot a\mbox{ \ for \ }x\in
\LL.
\end{equation}

An element $r$ in a Lie algebra $\LL$ is said to satisfy the
\textit{modern Yang-Baxter equation} (MYBE) if
\begin{equation}
\label{gvir-MYBE} x\cdot c(r)=0\mbox{ \ for all \ }x\in \LL.
\end{equation}

The main results of this paper is the following.
\begin{theo}\adddot
\label{gvir-main} {\rm(1)} Every Lie bialgebra structure on the Lie
algebra $\SS $ defined in $(\ref{gvir-equa-gvir-like})$ is a
triangular coboundary Lie bialgebra.

{\rm(2)} An element $r\in\SS $ satisfies CYBE in $(\ref{gvir-CYBE})$
if and only if it satisfies MYBE in $(\ref{gvir-MYBE})$.

{\rm(3)} Regarding $V=\SS \otimes\SS $ as an $\SS $-module under the
adjoint diagonal action of $\SS $ in $(\ref{gvir-diag})$, we have
$H^1(\SS ,V)=\Der(\SS ,V)/\Inn(\SS ,V)=0$.
\end{theo}
 \vskip10pt

\cl{\bf\S2. \ Proof of the main results}\setcounter{section}{2}
\setcounter{theo}{0}\setcounter{equation}{0} First we retrieve some
useful results from Drinfeld [D2], Michaelis [M1], Ng-Taft [NT] and
combine them into the following theorem.

\begin{theo}\adddot
\label{gvir-some} {\rm(1)} For a Lie algebra $\LL$ and
$r\in\Im(1-\tau)\subset \LL$, the tripple $(\LL,[\cdot,\cdot],
\D_r)$ is a Lie bialgebra if and only if $r$ satisfies MYBE
{\rm[D2]}.

{\rm(2)} Let $\LL$ be a Lie algebra containing two elements $a,b$
satisfying $[a, b] = b$, and set $r=a \otimes b -b \otimes a$. Then
$\D_{r}$ equips $\LL$ with the structure of a triangular coboundary
Lie bialgebra {\rm[M1]}.

{\rm(3)} For a Lie algebra $\LL$ and $r\in \Im(1-\tau)\subset \LL$,
we have {\rm[NT]}
\begin{equation}
\label{gvir-add-c} (1 + \xi + \xi^{2}) \cdot (1 \otimes \D) \cdot \D
(x) = x \cdot c (r) \mbox{ \ for all \ }x\in \LL.
\end{equation}
\end{theo}

We shall follow [SS, WS] closely to prove Theorem \ref{gvir-main}.

First Theorem \ref{gvir-main}(2) follows from the following more
general result.
\begin{lemm}\adddot
\label{gvir-lemma2} Denote by $\SS^{\otimes n}$ the tensor product
of $n$ copies of $\SS $. Regarding $\SS^{\otimes n}$ as an
$\SS$-module under the adjoint diagonal action of $\SS $, suppose
$c\in\SS^{\otimes n}$ satisfying $a\cdot c=0$ for all $a\in\SS $.
Then $c=0$.
\end{lemm}
\ni{\it Proof.} The lemma is obtained by using the same arguments in
the proof of [WS, Lemma 2.2]. \QED \vskip5pt Theorem
\ref{gvir-main}(3) follows from the following proposition.
\begin{prop}\adddot
\label{gvir-lemma3?} $\Der(\SS ,V)=\Inn(\SS ,V)$, where $V=\SS
\otimes\SS $.
\end{prop}
\ni{\it Proof.} We shall prove the result by several claims. Note
that $V=\oplus_{\a\in\G}V_\a$ is $\G$-graded with
$V_\a=\sum_{\b+\g=\a} \SS_\b\otimes\SS_\g$, where $\SS_\a=\F
L_\a\oplus \d_{\a,0}(\F\ptl_1+\F\ptl_2)$ for $\a\in\G$. A derivation
$D\in\Der(\SS ,V)$ is {\it homogeneous of degree $\a\in\G$} if
$D(V_\b) \subset V_{\a +\b}$ for all $\b \in\G$. Denote $\Der(\SS ,
V)_\a = \{D\in \Der(\SS , V) \,|\,{\rm deg\,}D =\a\}$ for $\a\in\G$.
\begin{clai}\adddot
\label{gvir-clai1} \rm Let $D\in\Der(\SS ,V)$. Then
\begin{equation}
\label{gvir-summable} \mbox{$D = \sum\limits_{\a \in\G} D_\a ,
\mbox{  \ where \ }D_\a \in \Der(\SS, V)_\a,$}
\end{equation}
which holds in the sense that for every $u \in\SS $, only finitely
many $D_\a (u)\neq 0,$ and $D(u) = \sum_{\a \in\G} D_\a(u)$ (we call
such a sum in (\ref{gvir-summable}) {\it summable}).
\end{clai}
\ \indent For $\a\in\G$, we define $D_\a$ as follows: For any
$u\in\SS _\b$ with $\b\in\G$, write $d(u)=\sum_{\g\in\G}v_\g\in V$
with $v_\g\in V_\g$, then we set $D_\a(u)=v_{\a+\b}$. Obviously
$D_\a\in\Der(\SS ,V)_\a$ and we have (\ref{gvir-summable}).
\begin{clai}\adddot
\label{gvir-clai2} \rm If $\a\ne0$, then $D_\a\in\Inn(\SS ,V)$.
\end{clai}
\ \indent Denote $T={\rm span}\{\ptl_1,\ptl_2\}$ and define the
nondegenerate bilinear map from $\G\times T\to\F$,
\begin{equation}\label{gvir-equa-nonde}
\ptl(\a)=\la\ptl,\a\ra=\la\a,\ptl\ra=a_1\a_1+a_2\a_2 \mbox{ for
}\a=(\a_1,\a_2)\in\G,\ \ptl=a_1\ptl_1+a_2\ptl_2\in T.
\end{equation}
By linear algebra, one can choose $\ptl\in T$ with $\ptl(\a)\not=0.$
Denote $a=(\ptl(\a))^{-1}D_{\a}(\ptl)\in \SS_{\a}.$ Then for any
$x\in \SS_{\b},\b \in \G,$ applying $D_{\a}$ to
$[\ptl,x]=\ptl(\b)x,$ using $D_{\a}(x)\in V_{\a+\b},$ We have
\begin{equation}\label{gvir-equa-add-1}
\ptl(\a+\b)D_{\a}(x) - x\cdot D_{\a}(\ptl)=\ptl \cdot D_{\a}(x)-
x\cdot D_{\a}(\ptl)=\ptl(\b)D_{\a}(x), \end{equation}
 i.e.,
$D_{\a}(x)=a_{\rm inn}(x).$ Thus $D_{\a}=a_{\rm inn}$ is inner.
%\hfill$\Box$

\begin{clai}\adddot
\label{gvir-clai3} \rm $D_0\in\Inn(\WW,V)$.
\end{clai}
\ \indent Choose a $\F$-basis
$\{\varepsilon_1,\varepsilon_2\}\subset \G$ of $\F.$ Define
$\ptl'_i\in T$ by $\langle \ptl'_i,\varepsilon_j \rangle =\d_{ij}.$
Let $\G'=\{(p,q)\in\F^2\,|\mbox{ }p\varepsilon_1+q\varepsilon_2\in
\G\}.$ Then $\Z^2\subset \G'.$ We write $L_{p,q}=L_{p\es_1+q\es_2}$,
and re-denote $\ptl'_i$ and $\G'$ by $\ptl$ and $\G$ respectively.
From (\ref{gvir-equa-gvir-like}), we have
$$
[L_{p,q},L_{p',q'}] = (qp'-pq')L_{p+q,p'+q'},\mbox{ \ }
 [\ptl_1,
L_{p,q}] = p L_{p,q},\mbox{ \ }  [\ptl_2,L_{p,q}] = q L_{p,q},
$$
for $(p,q),\,(p',q')\in\G\bs\{0\}$. The proof of this claim will be
done by several subclaims.
\begin{subc}\addbra \label{gvir-sub1} \rm $D_0(\ptl)=0$ for $\ptl\in T$.
\end{subc}
\par
To prove this, applying $D_0$ to $[\ptl,x]=\ptl(\b)x$ for
$x\in\SS_\b,\,\b\in \Gamma$, as in (\ref{gvir-equa-add-1}), we
obtain that \mbox{$x\cdot D_0(\ptl)=0$.} Thus by lemma 2.2,
$D_0(\ptl)=0$. \vskip4pt

\begin{subc}\addbra \label{gvir-sub2} \rm By replacing
$D_0$ by $D_0-u_{\rm inn}$ for some $u\in V_0$, we can suppose
$D_0(L_{p,q})=0$ for $p,q,p+q\in\{-1,0,1\}$.
\end{subc}
\vskip4pt
\par
We shall simplify notions by denoting
$$L^{p,q}_{r,s}=L_{p,q}\otimes
L_{r,s},\ \ L^{(i)}_{p,q}=\ptl_i\otimes L_{p,q},\ \
R^{(i)}_{p,q}=L_{p,q}\otimes \ptl_i \mbox { \ \ for \
}(p,q),(r,s)\in \G,\ i=1,2.$$ Denote by ${\rm Re}\, q$ the real part
of $q$ $\mbox{ for }q \in \F.$ Write
\begin{equation}
\label{gvir-equa2.5} D_0(L_{0,1}) =
\sum\limits_{p,q}c_{p,q}L^{p,q}_{-p,1-q} + c_1L_{0,1}^{(1)} +
d_1R_{0,1}^{(1)} + c_2L_{0,1}^{(2)} + d_2R_{0,1}^{(2)},
\end{equation}
 for some $c_{p,q},c_i, d_i\in
\F,$ where $\{ (p,q)\in\G\,|\mbox{ } c_{p,q}\ne 0\}$ is a finite
set. Note that
$$\begin{array}{lllllll} (L^{p,q-1}_{-p,1-q})_{\rm inn}(L_{0,1})&
= &
p(L^{p,q}_{-p,1-q} - L^{p,q-1}_{-p,2-q}),
\\[4pt]
(\ptl_2\otimes \ptl_2)_{\rm inn}(L_{0,1})& = & -R_{0,1}^{(2)}
-L_{0,1}^{(2)},\\[4pt]
(\ptl_1\otimes \ptl_2)_{\rm inn}(L_{0,1})& = & -L_{0,1}^{(1)},
\\[4pt]
(\ptl_2\otimes \ptl_1)_{\rm inn}(L_{0,1})& = & -R_{0,1}^{(1)}.
\end{array}$$
Using the above equations, by replacing $D_0$ by $D_0 - u_{\rm
inn}$, where $u$ is a combination of some $L^{p,q-1}_{-p,1-q},
\ptl_2\otimes \ptl_2, \ptl_1\otimes \ptl_2, \ptl_2\otimes \ptl_1,$
we can rewrite (\ref{gvir-equa2.5}) as (recall that $L_{0,0}=0$)
\begin{equation}
\label{gvir-equa2.6} D_0(L_{0,1}) = \sum\limits_{q\ne
0,1}c_{q}L^{0,q}_{0,1-q} + \sum\limits_{p\ne 0, \,0\leq {\rm Re}\,
q<1}c_{p,q}L^{p,q}_{-p,1-q} + cR_{0,1}^{(2)},
\end{equation}
 for some $c_{q},c_{p,q},
c \in \F,$ where $\{(0,q),(p,q)\in\G \,|\mbox{ } c_{q},c_{p,q}\ne
0\}$ is a finite set. Write
 \begin{equation}
 \label{gvir-equa2.7}
 D_0(L_{0,-1}) = \sum\limits_{q\ne
0,1}d_{q}L^{0,q-1}_{0,-q}
 +\sum\limits_{p\ne 0}d_{p,q}L^{p,q-1}_{-p,-q}
+ b_1 L_{0,-1}^{(1)} + f_1R_{0,-1}^{(1)} + b_2L_{0,-1}^{(2)} +
f_2R_{0,-1}^{(2)},
\end{equation}
 for some $d_q, d_{p,q}, b_i,
f_i\in \F,$ where $\{(0,q),(p,q)\in\G\,|\mbox{ } d_q,d_{p,q}\ne 0\}$
is a finite set. Applying $D_0$ to $[L_{0,1},L_{0,-1}] = 0,$ we have
$$
\mbox{$\sum\limits_{p\ne 0}$}d_{p,q}(pL^{p,q}_{-p,-q}-
pL^{p,q-1}_{-p,1-q})-b_2L^{0,1}_{0,-1} -f_2L^{0,-1}_{0,1}
=\mbox{$\sum\limits_{p\ne 0, \,0\leq {\rm
Re}\,q<1}$}c_{p,q}(pL^{p,q}_{-p,-q}-pL^{p,q-1}_{-p,1-q})
+cL^{0,1}_{0,-1}.
$$
Comparing the coefficients, we obtain $d_{p,q}=c_{p,q} \mbox{ for
}0\leq {\rm Re}\,q<1,\, p\ne 0, \mbox{ and }
 d_{p,q}=0 \mbox{ for  }p\ne 0,\, {\rm Re}\,q< 0
  \mbox{ or }{\rm Re}\,q\geq 1, \mbox{ and }b_2=-c, f_2=0.$
 Thus we can
rewrite (\ref{gvir-equa2.7}) as
\begin{equation}
\label{gvir-equa2.8} D_0(L_{0,-1}) = \sum\limits_{q\ne
0,1}d_{q}L^{0,q-1}_{0,-q}
 +\sum\limits_{p\ne 0, \,0\leq {\rm Re}\,q<1}c_{p,q}L^{p,q-1}_{-p,-q}
+ b_1 L_{0,-1}^{(1)} + f_1R_{0,-1}^{(1)} -cL_{0,-1}^{(2)}.
\end{equation}
% for some $d_q, c_{p,q},
% b_1, c, f_1
%\in \F,$ where $\{(0,q)\in\G~|\mbox{ } d_q\ne 0\}$ is a finite set.
Write
\begin{equation}
\label{gvir-equa2.9} D_0(L_{1,0}) =
\sum\limits_{p,q}e_{p,q}L^{p+1,q}_{-p,-q} + e_1 L_{1,0}^{(1)} +
e_2R_{1,0}^{(2)}+ e'_1R_{1,0}^{(1)} + e'_2L_{1,0}^{(2)},
\end{equation}
 for some $e_{p,q},e_i, e'_i\in \F,$ where $\{ (p,q)\in\G\,|\mbox{ }e_{p,q}\ne 0\}$ is a
finite set. Note that
$$\begin{array}{rcl}
(L^{0,p}_{0,-p})_{\rm inn}(L_{1,0})& = &
-p(L^{1,p}_{0,-p} - L^{0,p}_{1,-p}),\\[4pt]
(\ptl_1\otimes \ptl_1)_{\rm inn}(L_{1,0})& = & -R_{1,0}^{(1)} -
L_{1,0}^{(1)}.
\end{array}$$
Using these two equations, by replacing $D_0$ by $D_0 - u_{\rm
inn}$, where $u$ is a combination of some $L^{0,p}_{0,-p}{\ssc\,},\,
\ptl_1\otimes \ptl_1$ (this replacement does not affect the above
equations (\ref{gvir-equa2.6}), (\ref{gvir-equa2.8})), we can
rewrite (\ref{gvir-equa2.9}) as
\begin{equation}
\label{gvir-equa2.10} D_0(L_{1,0}) = \sum\limits_{p\ne
0}e_{p,q}L^{p+1,q}_{-p,-q} + e_1R_{1,0}^{(1)} + e_2R_{1,0}^{(2)} +
e'_2L_{1,0}^{(2)}.
\end{equation}
% for some $e_{p,q},e_1, e_2, e'_2\in \F,$ where
%$\{ e_{p,q}~|\mbox{ }e_{p,q}\ne 0\}$ is a finite set.
Applying $D_0$
to %the equation
$[L_{0,-1},[L_{0,1},L_{1,0}]]
%=[L_{0,-1},L_{1,1}]
=-L_{1,0},$ we \vs{-6pt}have
$$ \begin{array}{l}  \sum\limits_{p\ne
0}e_{p,q}\biggl( -(p\!+\!1)^2L^{p+1,q}_{-p,-q}\!+\!p(p\!+\!1)
L^{p+1,q+1}_{-p,-1-q}\!+\!p(p\!+\!1)L^{p+1,q-1}_{-p,1-q}
\!-\!p^2L^{p+1,q}_{-p,-q}\biggr) \! -\!e_1R_{1,0}^{(1)}
\! -\!e_2R_{1,0}^{(2)}  \\[12pt]
+e_2L^{1,1}_{0,-1}+e_2L^{1,-1}_{0,1}+e'_2L^{0,1}_{1,-1}
 +e'_2L^{0,-1}_{1,1}-e'_2 L_{1,0}^{(2)} -\sum\limits_{q\ne 0,1}c_{q}
 \biggl(qL^{1,q-1}_{0,1-q}-(q-1)L^{0,q}_{1,-q}\biggr)
-cR_{1,0}^{(2)}\\[12pt]
 -\! \sum\limits_{p\ne 0, \,0\leq {\rm
Re}\,q<1}\!\!c_{p,q}\biggl(
q(p\!+\!1)L^{p+1,q-1}_{-p,1-q}\!-\!pqL^{1+p,q}_{-p,-q}\!-\!p(q\!-\!1)L^{p,q-1}_{1-p,1-q}
\!+\!(q\!-\!1)(p\!-\!1)L^{p,q}_{1-p,-q}\biggr)\\[12pt]
 -\!\sum\limits_{q\ne
0,1}\!d_{q}\biggl(-(q\!-\!1)L^{1,q}_{0,-q}\!+\!qL^{0,q-1}_{1,1-q}\biggr)
\!-\!\sum\limits_{p\ne 0, \,0\leq {\rm Re}\,q<1}\!
c_{p,q}\biggl((p\!-\!q\!+\!1)L^{p+1,q}_{-p,-q}\!-(p\!-\!q)L^{p,q-1}_{1-p,1-q}\biggr)\\[12pt]
+b_1L^{1,1}_{0,-1}-b_1 L_{1,0}^{(1)}-f_1R_{1,0}^{(1)}+f_1
L^{0,-1}_{1,1} +cL_{1,0}^{(2)}
\\[12pt]
 =  -\sum\limits_{p\ne 0}e_{p,q}L^{p+1,q}_{-p,-q} -
e_1R_{1,0}^{(1)} - e_2R_{1,0}^{(2)} - e'_2 L_{1,0}^{(2)}\vs{-6pt}.
\end{array}$$
Comparing the coefficients of $R_{1,0}^{(1)}, \,L_{1,0}^{(1)},\,
R_{1,0}^{(2)}, \,L^{1,1}_{0,-1},\, L^{1,-1}_{0,1}, \,L^{0,1}_{1,-1},
\,L^{0,-1}_{1,1}$ respectively, we \vs{-6pt}obtain
\begin{equation}\label{gvir-equa2.11} f_1 = b_1 = c = 0,\mbox{ }~ e_2 = 2c_2 = 2d_{-1},\mbox{ }~ e'_2 = 2d_2 =
2c_{-1}\vs{-6pt}.\end{equation}
 Comparing the coefficients of
$L^{0,q}_{1,-q},\, L^{1,q}_{0,-q}$ with $q\ne 0,\pm1$ respectively,
we obtain $$(q-1)c_q = (q+1)d_{q+1},~~~~ (q+1)c_{q+1}=(q-1)d_q.$$
Note that $d_q=c_q=0 \mbox{ for } {\rm Re}\,q\gg 0 \mbox{ or }{\rm
Re}\,q\ll 0.$ \vs{-6pt}Thus the above equation forces
\begin{equation}\label{gvir-equa2.12} c_q = d_q = 0\ \mbox{ for
}q\ne 0,\vs{-6pt}1\end{equation}
 From (\ref{gvir-equa2.11}) and (\ref{gvir-equa2.12}), we have
$e'_2=0.$
 Comparing the coefficients
of $L^{p+1,q}_{-p,-q}$ with $p$ $\ne0,-1,\, {\rm Re}\,q<0 \mbox{ or
} {\rm Re}\,q\geq 1,$ we \vs{-6pt}obtain
\begin{equation}
\label{gvir-equa-add1-} e_{p,q-1}+e_{p,q+1}=2e_{p,q}\mbox{ \ \ for \
}p\ne0,-1,\, \ {\rm Re}\,q<0\mbox{ or }{\rm Re}\,q\geq 1\vs{-6pt}.
\end{equation}
%Setting $q$ to be $q+1,q+2,...,q+n$ if ${\rm Re\,}q\ge1$, or
%$q-1,q-2,...,q-n$ if ${\rm Re\,q}<0$, we solve from
Replacing $q$ by $q+n$ in (\ref{gvir-equa-add1-}) for $n\in\Z$, one
can \vs{-6pt}solve
$$ e_{p,q+n}=\biggl\{\begin{array}{ll} e_{p,q}+n(e_{p,q}-e_{p,q-1})&\mbox{if
\ }n\ge0\mbox{ and }{\rm Re\,}q\ge1,\\[4pt]
e_{p,q}-n(e_{p,q}-e_{p,q+1})&\mbox{if \ }n\le0\mbox{ and }{\rm
Re\,}q<0,
\end{array}\vs{-6pt}$$
for $p\ne0,-1$. However, $\{(p,q)\in\G\,|\,e_{p,q}\ne 0\}$ is a
finite set. We obtain $e_{p,q}=0\mbox{ for }p\ne 0,-1.$
 Comparing the coefficients of
$L^{p+1,q}_{-p,-q}$ with $p\ne0,-1,\, 0\leq {\rm Re}\,q<1$ , we
obtain $(p+1)c_{p,q}=pc_{p+1,q}.$  Thus $c_{p,q}=0 \mbox{ for
}p\ne0,\, 0\leq {\rm Re}\,q <1.$ Now (\ref{gvir-equa2.6}),
(\ref{gvir-equa2.8}) and (\ref{gvir-equa2.10}) \vs{-6pt}become
\begin{equation}\label{gvir-equa2.10+}
D_0(L_{0,1})=0, \mbox{ }~~ D_0(L_{0,-1})=0,\mbox{ }~~
D_0(L_{1,0})=\sum\limits_{q\ne0}e_qL^{0,q}_{1,-q}+eR_{1,0}^{(1)}\vs{-6pt},
\end{equation}
where $e=e_1,\,e_q=e_{-1,q}$.
%for some $e_q,e\in\F$ with
%$\{q\in\F\,|\,e_q\ne0\}$ being a finite set.
   \vs{-6pt}Write
\begin{equation}
\label{gvir-equa-d-0-l-1}
D_0(L_{-1,0})=\sum\limits_{p,q}f_{p,q}L^{p,q}_{-1-p,-q}+\tilde
f_1R_{-1,0}^{(1)}+\tilde f_2 L_{-1,0}^{(1)}+f'_1R_{-1,0}^{(2)}+
f'_2L_{-1,0}^{(2)},
\end{equation}
 for some $f_{p,q},\tilde
 f_i, f'_i\in
\F,$ where $\{(p,q)\in\G\,|\mbox{ }f_{p,q}\ne0\}$ is a finite set.
Applying $D_0$ to $[L_{0,-1},[L_{0,1},L_{-1,0}]]=-L_{-1,0}$, using
(\ref{gvir-equa2.10+}), we \vs{-7pt}obtain
$$\begin{array}{l}
  \sum\limits_{p,q}f_{p,q}
  \biggl(-p^2L^{p,q}_{-1-p,-q}+p(p+1)L^{p,q+1}_{-1-p,-1-q}+p(p+1)L^{p,q-1}_{-p-1,1-q}
-(p+1)^2L^{p,q}_{-1-p,-q}\biggr) \\[10pt]-\tilde f_1R_{-1,0}^{(1)}-\tilde  f_2
L_{-1,0}^{(1)}-f'_1R_{-1,0}^{(2)} -f'_1L^{-1,1}_{0,-1}-
f'_1L^{-1,-1}_{0,1} -f'_2 L^{0,1}_{-1,-1} -f'_2L^{0,-1}_{-1,1}- f'_2
L_{-1,0}^{(2)}\\[10pt]
 =  -\sum\limits_{p,q}f_{p,q}L^{p,q}_{-1-p,-q} -\tilde  f_1R_{-1,0}^{(1)}
  -\tilde f_2L_{-1,0}^{(1)}-
f'_1R_{-1,0}^{(2)} - f'_2L_{-1,0}^{(2)}\vs{-7pt}.
\end{array}$$
Comparing the coefficients of $L^{0,1}_{-1,-1}$, $L^{-1,-1}_{0,1}$
respectively, we obtain $f'_1=f'_2=0.$ Comparing the coefficients of
$L^{p,q}_{-1-p,-q}$ with $ p\ne 0,-1,$ we obtain
$f_{p,q-1}+f_{p,q+1}=2f_{p,q}.$ As in (\ref{gvir-equa-add1-}), by
noting that $\{(p,q)\in\G\,|\mbox{ }f_{p,q}\ne0\}$ is a finite set,
we obtain $f_{p,q}=0 \mbox{ for } p\ne0,-1.$ Now we can rewrite
(\ref{gvir-equa-d-0-l-1}) as
\begin{equation}\label{gvir-equa2.10++}
D_0(L_{-1,0})=\sum\limits_{q\ne0}f_{q}L^{0,q}_{-1,-q}
+\sum\limits_{q\ne0}f'_{q}L^{-1,q}_{0,-q}+fR_{-1,0}^{(1)}
+gL_{-1,0}^{(1)}, \end{equation} where
$f_{q}=f_{0,q},\,f'_{q}=f_{-1,q},\,f=\tilde f_1,\,g=\tilde f_2$.
Applying $D_0$ to $[L_{-1,0},L_{1,0}]$ $=0,$ we have
$$\begin{array}{ll}
&\sum\limits_{q\ne0}e_q(qL^{-1,q}_{1,-q}-qL^{0,q}_{0,-q})+eL^{1,0}_{-1,0}+fL^{-1,0}_{1,0}+gL^{1,0}_{-1,0}\\[12pt]
&=\sum\limits_{q\ne0}f_q(-qL^{1,q}_{-1,-q}+qL^{0,q}_{0,-q})+\sum\limits_{q\ne0}f'_q(-qL^{0,q}_{0,-q}
+qL^{-1,q}_{1,-q}).
\end{array}$$
Comparing the coefficients of $L^{-1,0}_{1,0},\,
L^{1,0}_{-1,0},\,L^{1,q}_{-1,-q},\,L^{-1,q}_{1,q}\mbox{ with }q\ne
0$ respectively, we obtain $f=f_q=0,\, e+g=0,\, e_q=f'_q.$ Thus
(\ref{gvir-equa2.10++}) \vs{-7pt}becomes
\begin{equation}\label{gvir-equa2.10+++}
D_0(L_{-1,0})=\sum\limits_{q\ne0}e_{q}L^{-1,q}_{0,-q}-eL_{-1,0}^{(1)}\vs{-7pt}.
\end{equation}
Applying $D_0$ to
$[L_{-1,0},[L_{1,0},L_{0,1}]]%=[L_{-1,0},-L_{1,1}]
=-L_{0,1},$ we \vs{-7pt}have
$$\begin{array}{ll}&\sum\limits_{q\ne0}e_{q}(qL^{-1,q}_{1,1-q}-(q-1)L^{0,q}_{0,1-q})
+eR_{0,1}^{(1)}
+eL^{1,1}_{-1,0}\\[12pt]&=\sum\limits_{q\ne0}e_{q}(-(q+1)L^{0,q+1}_{0,-q}+qL^{-1,q}_{1,1-q})
+eL_{0,1}^{(1)} +eL^{1,1}_{-1,0}\vs{-7pt}.
\end{array}$$
Comparing the coefficients of $R_{0,1}^{(1)},$
 we obtain $e=0.$ Comparing the coefficients of $L^{0,q}_{0,1-q}$ with
$q\ne 0,1,$ we obtain $(q-1)e_q=qe_{q-1}.$ Thus $e_q=0 $ for $q\ne
0$, and (\ref{gvir-equa2.10+}), (\ref{gvir-equa2.10+++})
\vs{-7pt}become
\begin{equation}
\label{gvir-equa2.14}
D_0(L_{1,0})=D_0(L_{-1,0})=D_0(L_{0,1})=D_0(L_{0,-1})=0\vs{-7pt}.
\end{equation}
From (\ref{gvir-equa2.14}) we can easily prove Subclaim 2).

\begin{subc}\addbra \label{gvir-sub3} \rm
$D_0(L_{p,q})=0$ for $(p,q)\in \Gamma\setminus\{0\}$.
\end{subc}
\vskip4pt
\par
Note that $L_{s,t}$ with $(s,t)\in  \Z^2\setminus \{0\}$ can be
generated by $\{L_{p,q}\,|\mbox{ }p,q,p+q\in\{0,\pm1\}\}$. From
Subclaim 2), we can easily deduct that $D_0(L_{p,q})=0$ for
$(p,q)\in \Z^2\setminus\{0\}.$
\par
For any element $(x,y)\in \G\setminus \Z^2,$ \vs{-7pt}write
\begin{equation}
\label{gvir-equa-??}
D_0(L_{x,y})=\sum\limits_{p,q}c_{p,q}L^{p,q}_{x-p,y-q}+a_1
L_{x,y}^{(1)}+b_1R_{x,y}^{(1)}+a_2
L_{x,y}^{(2)}+b_2R_{x,y}^{(2)}\vs{-7pt},
\end{equation}
for some $c_{p,q},a_i,b_i\in\F$.
 Applying $D_0$ to $[L_{0,-1},[L_{0,1},L_{x,y}]]=-x^2L_{x,y}$ and
comparing corresponding coefficients, we can obtain $c_{p,q}=0
\mbox{ for }p\ne 0,x.$ Similarly applying $D_0$ to
$[L_{-1,0},[L_{1,0},L_{x,y}]]=-y^2L_{x,y},$ we have $c_{p,q}=0
\mbox{ for }q\ne 0,y.$ Thus we can rewrite (\ref{gvir-equa-??}) as
$$D_0(L_{x,y})=eL^{x,0}_{0,y}+fL^{0,y}_{x,0}+a_1 L_{x,y}^{(1)}+b_1R_{x,y}^{(1)}
+a_2L_{x,y}^{(2)}+b_2R_{x,y}^{(2)}\mbox{ \ \ for \ }
e=c_{x,0},\,f=c_{0,y}.$$
 Applying
 $D_0$ to $[L_{-k,-1},[L_{k,1},L_{0,y}]]=-k^2y^2L_{0,y}\mbox { and }[L_{-1,-k},[L_{1,k},L_{x,0}]]=-k^2x^2L_{x,0}$ with
$k\gg 0,$ we can easily deduct that
 $D_0(L_{0,y})=D_0(L_{x,0})=0.$ Thus now we can assume that $xy\ne
 0.$ Applying $D_0$ to
 $[L_{-k,-1},[L_{k,1},L_{x,y}]]=-(x-ky)^2L_{x,y}$
 with $k\gg 0$ such that $ x-ky\ne 0,$ we have
$$\begin{array}{l}
-ex^2L^{x,0}_{0,y}+ekxyL^{x+k,1}_{-k,y-1}+ekxyL^{x-k,-1}_{k,y+1}-ek^2y^2L^{x,0}_{0,y}
-fk^2y^2L^{0,y}_{x,0}+fkxyL^{k,y+1}_{x-k,-1}\\[10pt]
+fkxyL^{-k,y-1}_{x+k,1} -fx^2L^{0,y}_{x,0}
+a_1k(x-ky)L^{k,1}_{x-k,y-1}+a_1k(x-ky)L^{-k,-1}_{x+k,y+1}-a_1(x-ky)^2L_{x,y}^{(1)}\\[10pt]
-b_1(x-ky)^2R_{x,y}^{(1)}+b_1k(x-ky)L^{x+k,y+1}_{-k,-1}+b_1k(x-ky)L^{x-k,y-1}_{k,1}+a_2(x-ky)L^{k,1}_{x-k,y-1}\\[10pt]
+a_2(x-ky)L^{-k,-1}_{x+k,y+1}-a_2(x-ky)^2L_{x,y}^{(2)}-b_2(x-ky)^2R_{x,y}^{(2)}
+b_2(x-ky)L^{a+k,y+1}_{-k,-1}\\[10pt]
+b_2(x-ky)L^{x-k,y-1}_{k,1}\\[10pt]=
-(x-ky)^2(eL^{x,0}_{0,y}+fL^{0,y}_{x,0}+a_1
L_{x,y}^{(1)}+b_1R_{x,y}^{(1)}+a_2 L_{x,y}^{(2)}+b_2R_{x,y}^{(2)}).
\end{array}$$
Comparing the coefficients of $L^{x,0}_{0,y}, L^{0,y}_{x,0},
L^{k,1}_{x-k,y-1}, L^{x+k,y+1}_{-k,-1}$ respectively, we have
\begin{equation}\label{gvir-equa2.16}
ekxy=fkxy=0,\mbox{ }~ a_1k+a_2=b_1k+b_2=0
\end{equation}
Since $xy\ne0$ and (\ref{gvir-equa2.16}) holds for all $k\in \Z$
with $k\gg 0$, we have $e=f=a_1=a_2=b_1=b_2=0.$ This proves Subclaim
3) and Claim \ref{gvir-clai3}. \begin{clai}\adddot
\label{gvir-clai4} \rm
 For every $D\in{\rm Der}(\SS,V)$, (\ref{gvir-summable}) is a finite sum.
 \end{clai}
\vskip4pt
\par By Claims \ref{gvir-clai2} and \ref{gvir-clai3}, we can suppose $D_\a=(v_\a)_{\rm
inn}$ for some $v_\a\in V_\a$ and $\a\in\G$. If
$\G'=\{\a\in\G\bs\{0\}\,|\,v_\a\ne0\}$ is an infinite set, by linear
algebra, there exists $\ptl\in T$ such that $\ptl(\a)\ne0$ for
$\a\in\G'$. Then $D(\ptl)=\sum_{\a\in\G'\cup\{0\}}\ptl\cdot v_\a=
\sum_{\a\in\G'}\ptl(\a)v_\a$ is an infinite sum, thus not an element
in $V$. This is a contradiction with the fact that $D$ is a
derivation from $\SS\to V$. This proves Claim \ref{gvir-clai4} and
the lemma.

\begin{lemm}\adddot \label{gvir-lemma3??} Suppose $r\in V$ such that
$a\cdot r\in {\rm Im}(1-\tau)$ for all $a\in \SS .$ Then $r\in {\rm
Im}(1-\tau)$
\end{lemm}

\ni{\it Proof.~}~(cf.~[SS], [WS]) First note that $\SS \cdot {\rm
Im}(1-\tau)\subset  {\rm Im}(1-\tau).$  We shall prove that after a
number of steps in each of which $r$ is replaced by $r - u$ for some
$u\in {\rm Im}(1-\tau),$ the zero element is obtained and thus
proving that $r\in{\rm Im}(1-\tau)$.
 Write $r=\sum_{x\in \Gamma}r_x.$ Obviously,
 \begin{equation}
 \label{gvir-equa2.17}
 r\in {\rm Im}(1-\tau)\ \ \Longleftrightarrow \ \ r_x \in {\rm Im}(1-\tau)
 \mbox{ for all } x\in \Gamma.
 \end{equation}
 For any $x'\ne0$, choose $\ptl\in T$ such that $\ptl(x')\ne 0$.
 Then $\sum_{x\in \Gamma}\ptl(x)r_x=\ptl\cdot r\in {\rm
 Im}(1-\tau).$ By (\ref{gvir-equa2.17}), $\ptl(x)r_x\in {\rm Im}(1-\tau),$ in
 particular, $r_{x'}\in {\rm Im}(1-\tau).$ Thus by replacing $r$
 by $r-\sum_{0\ne x\in \Gamma}r_x,$ we can suppose
 $r=r_0\in V_0.$
 Now we can write
 \begin{equation}
 r=\sum_{p,q}c_{p,q}L^{p,q}_{-p,-q}+c_1\ptl_1\otimes \ptl_1+c'_{1}\ptl_1\otimes
 \ptl_2+c_2\ptl_2\otimes \ptl_1+c'_2\ptl_2\otimes \ptl_2,
 \end{equation}
 for some $c_{p,q},c_i,c'_i\in\F$.
Choose any total order on $\F$ compatible with its group structure.
 Since $v_{p,q}:=L^{p,q}_{-p,-q}-L^{-p,-q}_{p,q}\in {\rm Im}(1-\tau),$ by replacing $r$ by $r-u$, where
 $u$ is a combination of some $v_{p,q}.$ We can suppose
 \begin{equation}\label{gvir-equa2.19}
 c_{p,q}\ne 0\ \ \Longrightarrow\ \  p>0 \mbox{ or } p=0, q>0.
 \end{equation}
 First assume that $c_{p,q}\ne0$ for some $p,q.$
Choose $s, t>0$ such that $sq-p{\ssc\,}t\ne 0.$ Then we see that the
term $L^{p+s,q+t}_{-p,-q}$ appears in $L_{s,t}\cdot r,$ but
(\ref{gvir-equa2.19}) implies that the term $L^{-p,-q}_{p+s,q+t}$
does not appear in $L_{s,t}\cdot r,$  a contradiction with the fact
that $L_{s,t}\cdot r\in {\rm Im}(1-\tau).$ Now write
$r=c_1\ptl_1\otimes \ptl_1+c'_{1}\ptl_1\otimes
 \ptl_2+c_2\ptl_2\otimes \ptl_1+c'_2\ptl_2\otimes \ptl_2.$
  Then from
  $$
  \begin{array}{rcl}
  L_{1,0}\cdot r& = &-c_1R_{1,0}^{(1)}-c_1
  L_{1,0}^{(1)}-c'_1 R_{1,0}^{(2)}-c_2
  L_{1,0}^{(2)}~~\in
  ~{\rm Im}(1-\tau),\\[12pt]
  L_{0,1}\cdot r& = &-c'_1L_{0,1}^{(1)}
  -c_2 R_{0,1}^{(1)}-c'_2 R_{0,1}^{(2)}-c'_2L_{0,1}^{(2)}~~\in ~{\rm
  Im}(1-\tau),
  \end{array}
$$ we obtain that $c_1=0,\, c'_1+c_2=0,\, c_2=0.$ Thus $r\in {\rm
Im}(1-\tau).$ This proves the lemma.

\vskip10pt \ni{\it Proof of Theorem \ref{gvir-main}(1).} Let $(\SS
,[\cdot,\cdot],\D)$ be a Lie bialgebra structure on $\SS $. By
(\ref{gvir-bLie-d}), (\ref{gvir-deriv}) and Theorem
\ref{gvir-main}(3), $\D=\D_r$ is defined by (\ref{gvir-D-r}) for
some $r\in\SS \otimes\SS $. By (\ref{gvir-cLie-s-s}), ${\rm
Im}\,\D\subset{\rm Im}(1-\tau)$. Thus by Lemma \ref{gvir-lemma3??},
$r\in{\rm Im}(1-\tau)$.
%, namely, $r=\sum_{i=1}^m (L_i\otimes b_i-b_i\otimes L_i)$ for some $L_i,b_i\in\WW$.
Then by (\ref{gvir-cLie-j-i}), (\ref{gvir-add-c}) and Theorem
\ref{gvir-main}(2) show that $c(r)=0$. Thus Definition
\ref{gvir-def2} says that $(\SS ,[\cdot,\cdot],\D)$ is a triangular
coboundary Lie bialgebra. \QED

\vskip12pt

\cl{\bf References}\vs{0pt}

\vskip3pt\small
\parindent=8ex\parskip=0pt\baselineskip=0pt

\re{D1} V.G. Drinfeld, Hamlitonian structures on Lie group,Lie
algebras and the geometric meaning of classical Yang-Baxter
equations, {\it Soviet Math. Dokl.} {\bf27}(1) (1983), 68-71.

\re{D2} V.G. Drinfeld, Quantum groups, in: {\it Proceeding of the
International Congress of Mathematicians}, Vol.~1, 2, Berkeley,
Calif.~1986, Amer.~Math.~Soc., Providence, RI, 1987, pp.~798-820.

\re{DZ} D.~Dokovic, K.~Zhao, Derivations, isomorphisms and
  second cohomology of generalized Block algebras,
  {\it Algebra Colloquium} {\bf3} (1996), 245--272.

%\re{DZ1} D. Dokovic, K. Zhao, Generalized Cartan type
%$S$ Lie algebras in characteristic zero, {\it J.
%Algebra} {\bf193} (1997), 643--664.

%\re{F} R. Farnsteiner, Derivations and central extensions of
%finitely generalized Lie algebras, {\it J. Algebra} {\bf 118}
%(1988), 33--45.

\re{LT} W.~Lin, S.~Tan, Nonzero level Harish-Chandra modules over
the Virasoro-like algebra, {\it J.~Pure Appl.~Algebra} (2006), in
press.

\re{MJ} D.~Meng, C.~Jiang, The derivation algebra and the universal
central extension of the $q$-analog of the Virasoro-like algebra,
{\it Comm. Algebra} {\bf 26} (1998), 1335--1346.

 \re{M1} W. Michaelis, A
class of infinite-dimensional Lie bialgebras containing the Virasoro
algebras, {\it Adv. Math.} {\bf107} (1994), 365--392.

\re{M2} W. Michaelis, The dual Poincare-Birkhoff-Witt theorem, {\it
Adv. Math.} {\bf57} (1985), 93--162.

%\re{M3} W. Michaelis, Lie coalgebras, {\it Adv. Math.} {\bf38}
%(1980), 1--54.

\re{NT} S.-H. Ng, E.J. Taft, Classification of the Lie
bialgebra structures on the Witt and Virasoro algebras,
{\it J. Pure Appl.~Algebra} {\bf151} (2000), 67--88.

\re{N} W.D. Nichols, The structure of the dual Lie coalgebra
of the Witt algebra, {\it J. Pure Appl. Algebra} {\bf 68} (1990),
395--364.

%\re{S}
%G.~Song, The structure if infinite dimensional non-graded Lie algebras and
%Lie superalgebras of $W$-type and the related problems, {\it Ph.~D.~Thesis,}
%Shanghai Jiaotong University (2005).

\re{SS} G. Song, Y. Su, Lie bialgebras of generalized Witt type,
{\it Science in China A}, in press (see arXiv:math.QA/0504168).

\re{SX} Y. Su, X. Xu, Structure of divergence-free Lie algebras,
{\it J. Algebra} {\bf 243} (2001), 557--595.

\re{T} E.J. Taft, Witt and Virasoro algebras as Lie bialgebras, {\it
J. Pure Appl.~Algebra} {\bf87} (1993), 301--312.

\re{WS} Y. Wu, Y. Su, Nongraded Lie bialgebras of generalized Witt
type, to appear.

\re{X1} X. Xu, New generalized simple Lie algebras of Cartan type
over a field with characteristic 0,  {\it J.~Algebra} {\bf 224}
(2000), 23--58.

\re{X2} X.~Xu, Novikov-Poisson algebras,
  {\it J.~Algebra} {\bf190} (1997), 253--279.

\re{X3} X.~Xu, Generalizations of Block algebras,
  {\it Manuscripta Math.} {\bf100} (1999), 489--518.

\re{ZZ} H. Zhang, K. Zhao, Representations of the Virasoro-like
algebra and its q-analog, {\it Comm. Algebra} {\bf24} (1996),
4361--4372.

\re{Z1}  K.~Zhao, A class of infinite dimensional simple Lie
algebras,  {\it J.~London Math.~Soc. (2)}, {\bf62} (2000), 71--84.

\re{Z2} K. Zhao, Generalized Cartan type $S$ Lie algebras in
characteristic zero, II, {\it Pacific J. Math.} {\bf 192} (2000),
431--454.

\re{ZM} L.~Zhu, D.~Meng, Structure of degenerate Block algebras,
  {\it Algebra Colloquium} {\bf10} (2003), 53--62.

\end{document}